\documentclass{amsart}%
\font \Bbbten=msbm10 \font \Bbbsev=msbm7
 \font \Bbbfiv=msbm5
\newfam\Bbbfam \def \Bbb {\fam\Bbbfam\Bbbten}\textfont\Bbbfam =\Bbbten
\scriptfont\Bbbfam=\Bbbsev \scriptscriptfont\Bbbfam=\Bbbfiv

\font\de=cmssi10

\begin{document}


\newcounter{theorem}[section]
\newtheorem{defi}[theorem]{\sc Definition}
\newtheorem{lema}[theorem]{\sc Lemma}
\newtheorem{prop}[theorem]{\sc Proposition}
\newtheorem{cor}[theorem]{\sc Corrollary}
\newtheorem{teo}[theorem]{\sc Theorem}
\newtheorem{obs}[theorem]{\sc Remark}
\newtheorem{ques}[theorem]{\sc Question}
\newtheorem{question}{\sc Question}
\def\bp{\noindent{\it Proof. }}
\def\ep{\noindent{\hfill $\fbox{\,}$}\medskip\newline}
\renewcommand{\theequation}{\arabic{section}.\arabic{equation}}
\renewcommand{\thetheorem}{\arabic{section}.\arabic{theorem}}
\newcommand{\eps}{\varepsilon}
\newcommand{\disp}[1]{\displaystyle{\mathstrut#1}}
\newcommand{\fra}[2]{\displaystyle\frac{\mathstrut#1}{\mathstrut#2}}
\newcommand{\Z}{\Bbb Z}
\newcommand{\N}{\Bbb N}
\def\to{\mathop{\rightarrow}}
\def\ord{\mathop{\rm ord}}
\def\diam{\mathop{\rm diam}}
\def\cc{\mathop{\rm cc}}
\title{Continuum-wise expansive homeomorphisms on Peano continua}
\author{Jana Rodriguez Hertz}
\date{}
\maketitle
\begin{abstract}
On a Peano continuum, all local stable and unstable components of
a continuum-wise expansive homeomorphism are non trivial. In
particular, there is sensitive dependence on initial conditions.
 This generalizes results in \cite{h,l} about
lack of Lyapunov stable points (weak sinks) and existence of non
trivial stable and unstable components for expansive
homeomorphisms on Peano continua.
\par We also use this fact to generalize a result in
\cite{ktt}: a Peano curve $X$ admitting a continuum-wise expansive
homeomorphism is nowhere rim-countable. However, it is not
resolved the question of whether such a dynamics could be possible
in a locally planar Peano curve. Some other questions are posed.
\end{abstract}
 A homeomorphism $f$ of a compact metric space $(X,d)$ is {\de continuum-wise expansive (cw-expansive)}
 if there exists a constant $\alpha>0$ such that $\sup_{n\in\Z}{\rm diam}\: f^n(C)>\alpha$ for all non trivial continua $C$. This notion was introduced in \cite{kato_cw}, and is a
generalization of {\de expansive} homeomorphisms: those verifying,
for some fixed $\alpha>0$, that
$\sup_{n\in\Z}d(f^n(x),f^n(y))>\alpha$ for all $x\ne y$. Both
notions are topological invariants. Examples of cw-expansive
homeomorphisms that are not expansive are non-expansive homeomorphisms of totally disconnected spaces.\par %
We shall use the following notations: for $\eps>0$ we define the
{\de local stable and unstable set of $x$}, respectively, as:
$$W^s_\eps(x)=\{y\in X:\: \sup_{n\in{\Bbb N}_0}d(f^n(x), f^n(y))\leq\eps\};\quad W^u_\eps(x)=\{y\in X:\: \sup_{n\in{\Bbb N}_0}d(f^{-n}(x), f^{-n}(y))\leq\eps\}$$
These sets are important in analyzing and classifying expansive
dynamics, though in the case of cw-expansive homeomorphisms, it
seems more reasonable to study the {\de stable and unstable
component of $x$}, that is:
$$CW^\sigma_\eps(x)=cc(W^\sigma_\eps(x),x)\qquad \sigma=s,u$$
We obtain the following result for cw-expansive homeomorphisms
acting on Peano continua:
\begin{teo} \label{piezas.conexas.cw} If $f$ is a cw-expansive
homeomorphism on a Peano continuum $X$, then for each $\eps>0$
there exists $\delta>0$ such that $\inf_{x\in X}{\rm
diam}\:CW^\sigma_\eps(x)\geq\delta$ for both $\sigma=s,u$\end{teo}
We recall that a Peano continuum is a locally connected continuum.
In particular, we have:
\begin{teo}\label{no.sinks.cw} Cw-expansive homeomorphisms on
Peano continua have sensitive dependence on initial conditions.
\end{teo}
In \S\ref{condiciones generales}, sufficient conditions are found
for a point in a general compactum to verify
$CW^s_\eps(x)\ne\{x\}$.
 No example of a cw-expansive
homeomorphisms of a Peano curve is known; and, in fact it is an
open question whether there is a Peano curve admitting an
expansive homeomorphism. In \cite{ktt} it was shown, by using
techniques in \cite{h,l} that a Peano curve admitting an expansive
homeomorphism is nowhere planar, nowhere rim-countable and
consists of a finite number of cyclic elements. We obtain here a
cw-expansive partial version of this:
\begin{teo}\label{nowhere.planar}A Peano curve admitting a
cw-expansive homeomorphism is nowhere rim countable. As a
consequence, it consists of a finite number of cyclic elements.
Also, there are at
most countably many local cut points, and every open set contains a simple close curve.\end{teo}%
We must remark, however, that the following question remains open
\begin{ques}
Is a Peano curve admitting a cw-expansive homeomorphism nowhere
locally planar?
\end{ques}
Observe that obtaining non trivial stable and unstable components
is not enough to apply techniques in \cite{ktt}, because in
cw-expansive dynamics there is not much information about the sets
$CW^s_\eps(x)\cap CW^u_\eps(y)$, except that it is a totally
disconnected set, and hence no local product
structure should be expected a priori.\par%
We also state the following questions:
\begin{ques}
Is $CW^\sigma_\eps(x)$ locally connected for cw-expansive
homeomorphisms acting on a Peano continua?
\end{ques}
This question has not been answered yet, even for expansive
homeomorphisms on manifolds.
\begin{ques}
Let $X$ be a Peano continuum and $f$ an expansive (cw-expansive) homeomorphism.
Is $\dim_xX$ constant over $X$? Maybe it can help considering a
transitive $f$.
\end{ques}
\begin{ques}\label{pregunta.sink.en.cuenca de periodico}
Let $X$ be a continuum and $f$ an expansive (cw-expansive)
homeomorphism. Does every sink belong to some $W^s(p)$, with $p$
periodic? (Look at \S\ref{basic definitions} for a definition of
sink) The point $p$ is not necessarily a sink, see Example in
\S\ref{ejemplo.1}.
\end{ques}
\begin{ques}\label{pregunta conexion local pozos}
Let $X$ be a continuum, and $f$ an expansive homeomorphism. Are
sinks and sources points of local connectedness of $X$?
\end{ques}
 {\em Acknowledgements. } I thank Ed Tymchatyn for kindly answering some questions.


\section{Stable and unstable components}\label{basic definitions}\label{conexion+estabilidad}
A point in a compact metric space $X$ is called a {\de weak sink}
or, equivalently, a Lyapunov stable point if $W^s_\eps(x)$ is a
neighborhood of $x$ for each $\eps>0$. A weak sink $x$ is a {\de
sink} if, moreover, $W^s_\eps(x)\subset W^s(x)$, that is, if all
points $y\in W^s_\eps(x)$ verify that $d(f^n(x),f^n(y))\to 0$ when
$n\to \infty$. Sources and weak sources are defined analogously.
Let us mention that in the literature it is usually required that
sinks be periodic, this notion is more general.\par%
Observe that, in an expansive setting, weak sinks are sinks, and
this fact still holds for cw-expansive homeomorphisms on Peano
continua. A map $f$ has {\de sensitive dependence on initial
conditions} if there are no weak sinks on $X$.\par%

In this section Theorem \ref{piezas.conexas.cw} is proved. It is a
standard fact that a non-wandering sink is periodic. Indeed, the
proof follows as in an expansive setting for any compact $X$. We
shall see that, moreover, weak sinks are always sinks and
non-wandering in a Peano continuum. This, together with the fact
that sinks form an open set, ends the proof of Theorem \ref{no.sinks.cw}.\par%
Theorem \ref{piezas.conexas.cw} follows by considering the
connected component of the set $W^s_{\eps,n}(x)$ (see below) that
is mapped into $B_\delta(f^n(x))$ by $f^n$. This is a connected
neighborhood of $x$ that must meet the boundary
of $W^s_{\eps,n}(x)$, for otherwise a source would appear in $\omega(x)$.\par%
Proposition below states that each continuum of diameter $\eps>0$
becomes $\alpha$-distinguishable in at most $N_\eps$ forward or
backward iterates, for some constant $N_\eps>0$ depending on $f$.
This statement follows from compactness of $X$.
\begin{prop}\label{cw-uniform} For all
$\eps>0$ there exists $N_\eps>0$ such that $\disp{\max_{|n|\leq
N_\eps} \diam f^n(C)>\alpha}$ for any continuum $C$ satisfying
$\diam C\geq\eps$.
\end{prop}
Let us define:
$$W^s_{\eps,n}(x)=\{y\in X: d(f^r(x), f^r(y))\leq\eps\quad
r=0,\dots,n\}\quad \mbox{and}\quad
CW^s_{\eps,n}(x)=cc(W^s_{\eps,n}(x),x)$$%
We have, as an immediate corollary of Proposition
\ref{cw-uniform}:
\begin{cor}\label{corolario.cw.eps} For all $x\in X$ and
$0<\eps<\alpha/2$ we have, if $n\geq N_\eps$, that:
\begin{enumerate}
\item \label{item1}$f^n(CW^s_{\alpha/2}(x))\subset
CW^s_\eps(f^n(x))$%
\item \label{item2} $CW^s_\eps(x)\subset W^s(x)$%
\item \label{item3} $CW^s_\eps(x)=\cc(W_{\alpha/2}^s(x)\cap
W_{\eps,n}^s(x),x)$.%
\item\label{item4} $x$ is a sink $\Leftrightarrow$
$CW_{\alpha/2}^\sigma(x)$ is a neighborhood of $x$.
\end{enumerate}
\end{cor}
Indeed, $C_n=f^n(CW^s_{\alpha/2}(x))$ is a continuum satisfying
$\diam f^k(C_n)\leq \alpha$ for all $|k|\leq N_\eps$, hence $\diam
C_n<\eps$. To see the less trivial inclusion $\supset$ in Item
\ref{item3}., take $D_n=\cc(W_{\alpha/2}^s(x)\cap
W_{\eps,n}^s(x),x)$. Then on one hand $D_n\subset
CW_{\eps,n}^s(x)$, while, in view of Item \ref{item1},
$f^n(D_n)\subset CW^s_\eps(f^n(x))$, so $D_n\subset
CW^s_\eps(x)$.\ep%
Observe that all statements in Corollary \ref{corolario.cw.eps}
hold in a general compactum $X$, except, in principle,
$\Rightarrow$ in Item \ref{item4}. See Question \ref{pregunta
conexion local pozos}.
\begin{prop} \label{conjunto.weak.sinks.cw} Let $X$ be a Peano
continuum and $f$ a cw-expansive homeomorphism. Then the set
$\mathcal{S}$ of sinks is an open set, and $\mathcal{S}\cap
\Omega(f)$ consists of periodic points.
\end{prop}
\bp Let $x$ be a sink, and take a connected open neighborhood of
$x$, $U\subset CW_{\alpha/4}^s(x)$. Then $U$ consists of sinks.
Indeed, triangular inequality implies $U\subset
CW^s_{\alpha/2}(y)$ for all $y\in U$, hence from item \ref{item4}
in Corollary above, we get $U$ consists of sinks. This shows the
sets of sinks is open.\par%
Assume now that $x$ is non-wandering, and take $M>N_{\alpha/4}$ so
that $f^M(U)\cap U\ne\emptyset$. But $f^M(CW^s_{\alpha/2}(x))$ is
contained, by Corollary above, in $CW^s_{\alpha/4}(f^M(x))$, so
choice of $M$ implies that $f^M(CW^s_{\alpha/2}(x))$ is in fact
contained in $CW^s_{\alpha/2}(x)$. Now, as $\diam
f^{kM}(CW^s_{\alpha/2}(x))\to 0$, we have
$\bigcap_kf^{kM}(CW^s_{\alpha/2}(x))$ must consist of only one
(periodic) point. But there is at most one non-wandering point in
$CW^s_{\alpha/2}(x)$, hence $x$ is periodic. \ep
The non existence of sinks will be proved by showing that all
sinks are, in fact, non-wandering in the presence of local
connectedness. We see that the basins of attraction of all sinks
$f^{-n}(x)$ contain a ball of uniform radius. Being $f^{-n}(x)$
infinitely many would imply they ``collapse", hence $x$ would be,
in fact, periodic by arguments in previous Proposition. Indeed,
basins of attractions cannot shrink too much in the past, or else
one would obtain a continuum not satisfying Proposition
\ref{cw-uniform}.
\begin{prop}\label{pozo no errante}
In a Peano continuum $X$, all sinks are non-wandering, hence
periodic.
\end{prop}
\bp There exist $\rho>0$ and $N\geq 0$ verifying
$\cc(B_\rho(f^{-n}(x)), f^{-n}(x))\subset CW_{\alpha/2}^s
(f^{-n}(x))$ for all $n\geq N$. Indeed, let $\eps>0$ be such that
$B_\eps(x)\subset CW_{\alpha/2}^s(x)$, and consider $\rho>0$ so
that $B_\rho(y)\subset W^s_{\eps,N_\eps}(y)$ for all $y\in X$.
Observe that $W_{\eps,n}^s(f^{-n}(x))\subset
CW_{\alpha/2}^s(f^{-n}(x))$ for all $n\geq0$. If we had
$\cc(B_\rho(f^{-n}(x)),f^{-n}(x))\not\subset
W_{\eps,n}^s(f^{-n}(x))$ for some $n\geq 2N_\eps$, then the
continuum $D=\cc(B_\rho(f^{-n}(x))\cap
W_{\eps,n}^s(f^{-n}(x)),f^{-n}(x))$ would satisfy $D\cap \partial
W_{\eps,n}^s(f^{-n}(x))\ne\emptyset$, hence $\diam f^k(D)= \eps$
for some $k\geq N_\eps$. This would contradict Proposition
\ref{cw-uniform}, since $\diam f^{j+k}(D)\leq\alpha$ for all
$|j|\leq N_\eps$. So $\cc(B_\rho(f^{-n}(x)),f^{-n})\subset
W_{\eps,n}^s(f^{-n}(x))$ for all $n\geq N_\eps$, and the claim
follows.\par%
Given $\gamma>0$, take $n>m>N_\gamma$ so that $f^{-m}(x)\in
\cc(B_\rho(f^{-n}(x)),f^{-n}(x))$. Corollary
\ref{corolario.cw.eps}.\ref{item1} implies
$f^{n-m}(x)\in CW^s_\gamma(x)$, and hence $x\in\omega(x)$.\ep%
In case $X$ is a Peano continuum we obtain, by putting together
Propositions \ref{conjunto.weak.sinks.cw} and \ref{pozo no
errante}, that, were the set of sinks not empty, there would be an
open set of periodic points. In particular, there would be two
different periodic points in the same basin of attraction. This
proves Theorem
\ref{no.sinks.cw}.\par%
Non triviality of stable and unstable components follows now from
this fact. However, sensitivity to initial conditions is not
enough, in a general continuum, to guarantee that these components
are nontrivial. Example in \S\ref{ejemplo.1} illustrates this.\par%
In order to prove Theorem \ref{piezas.conexas.cw}, it will be seen that
 $CW_{\eps,n}^s(x)\cap\partial B_\delta(x)\ne \emptyset$ for all
$n\geq K$. The proof follows then from the fact that
$CW^s_{\eps,n}(x)$ is a decreasing sequence converging to
$CW_\eps^s(x)$.\par%
Let $D_n=\cc(W^s_{\eps,n}(x)\cap f^{-n}(B_\delta(f^n(x))),x)$,
where $\delta>0$ verifies $B_\delta(y)\subset
W^s_{\eps,N_\eps}(y)\cap W^u_{\eps,N_\eps}(y)$ for all $y\in X$.
Observe that $D_n\cap \partial W^s_{\eps,n}(x)$ is not empty for
all sufficiently large $n>0$. Indeed, there exists $K\geq 2N_\eps$
such that $\cc(B_\delta(f^n(x)),f^n(x))\not\subset
f^n(W^s_{\eps,n}(x))=W^u_{\eps,n}(f^n(x))$ for all $n\geq K$, for
otherwise $\omega(x)$ would consist of sources.\par%
Hence for all $n\geq K$ there exists some $0\leq k\leq n$ such
that $\diam f^k(D_n)=\eps$. From this it follows that $k<N_\eps$
since, by choice of $\delta$, $k\leq n-N_\eps$; and by Proposition
\ref{cw-uniform} $k\notin [N_\eps, n-N_\eps ]$. Choice of $\delta
$ again, implies $D_n\cap \partial B_\delta(x)$ is not empty.\ep
\subsection{Sufficient conditions in general compact sets}\label{condiciones generales}
The following propositions are valid for a cw-expansive
homeomorphism $f$ of a compact metric space $X$.
\begin{prop}\label{alfa.pozos}
If $x$ is a sink such that $\alpha(x)$ contains a point of local
connectedness, then $x$ is periodic.
\end{prop}
\bp Let $z=\lim_kf^{-n_k}(x)\in \alpha(x)$ be a point of local
connectedness in $\alpha(x)$, and let us find $\rho>0$ such that
$\cc(B_\rho(z),z)\subset CW_{\alpha/2}^s (f^{-n_k}(x))$ for all
$k$.\par%
The proof follows as in Proposition \ref{pozo no errante},
considering, for all $f^{-n_k}(x)\in \cc(B_\rho(z),z)$, the
continuum $D=\cc(B_\rho(f^{-n_k})\cap
W_{\eps,n}^s(f^{-n_k}(x)),f^{-n_k}(x))$.\ep%
Related to next Proposition is the second example in
\S\ref{ejemplo.1}, which shows there can be an open set of points
verifying $W^s_\eps(x)=W^u_\eps(x)=\{x\}$, even if $X$ is a
continuum.
\begin{prop} For each $\eps>0$ there exists $\delta>0$ such that
if $\omega(x)$ contains a point of local connectedness, then
either $x$ is a periodic source or $\diam CW^s_\eps(x)\geq\delta$.
\end{prop}
Let $z=\lim_k f^{n_k}(x)$ be a point of local connectedness of
$\omega(x)$ and let $\delta>0$ be as in the proof of Theorem
\ref{piezas.conexas.cw}. If there is a subsequence such that
$CB_\delta(f^{n_k}(x))\not\subset W_{\eps,n_k}^s(f^{n_k}(x))$,
then $CW^s_\eps(x)$ is not trivial. The proof is as in Theorem
\ref{piezas.conexas.cw}. Moreover, $\diam
CW^s_\eps(x)>\delta$.\par%
Otherwise, $z$ is a source verifying  $CB_\delta(z)\subset
W^u_\eps(z)$. Argument in Proposition \ref{pozo no errante} shows
$z$ is periodic. Since $z\in\omega(x)$ there exists $f^n(x)\in
 CW^s_\eps(z)$. Being $z$ a source implies $f^n(x)=z$.\ep
\section{Cw-expansive homeomorphisms on Peano curves}
I thank R. Ures the suggestion of the following:
\begin{prop}
If $X$ is a Peano continuum admitting a cw-expansive
homeomorphism, then for each $x\in X$: $\ord(x)=2^{\aleph_0}$.
That is, $X$ is nowhere rim-countable.
\end{prop}
\bp Let $x\in X$, $0<\eps<\alpha/4$ and consider $U$ a connected
neighborhood of $x$ such that $U\subset B_\delta(x)$, where
$0<\delta<\frac14 \inf_{w\in X} \diam CW^u_\eps(w)$. Call
$s=\cc(W^s_\eps(x)\cap U,x)$, and consider, for all $z\in\partial
U$ the closed set $B_z=s\cap CW^u_\eps(z)$. $B_z$ has empty
interior in the topology of $s$, due to cw-expansiveness. But $s=
\bigcup_{z\in\partial U}B_z$. If $\partial U$ were countable then
$s$ would be a countable union of closed sets with empty interior,
which is absurd.\ep As a consequence of this, we can also
conclude:
\begin{cor}
If $X$ is a Peano curve admitting a cw-expansive homeomorphism,
then (1) All points in $X$ belong to a continuum of convergence.
(2) Every open set of $X$ contains a simple closed curve. (3) The
set of local cut points is at most countable. (4) $X$ has a finite
number of cyclic components.
\end{cor}
(1) For any $x_0\in X$ and $\eps>0$, there is a sequence $x_n\in
CW^s_\eps(x_0)$, with $x_n\to x_0$ such that
$\{CW^u_\eps(x_n)\}_{n\geq0}$ is a pairwise disjoint collection.
The Hausdorff limit of $CW^u_\eps(x_n)$ is a non trivial continuum
containing $x_0$. (2) The existence of a neighborhood basis
$\{K_\lambda(x)\}$ of a point $x$, such that each $K_\lambda(x)$
is a dendrite, would imply $x$ is a regular point, so each open
set contains a close curve. Item (3) is implied by results in
\cite{w}. Item (4) is shown as in \cite{ktt}.\ep
\section{Example}\label{ejemplo.1}
This is an example of an expansive dynamics in a continuum $X$,
with sensitivity to the initial conditions, but not verifying the
result stated in Theorem \ref{piezas.conexas.cw}. In this example,
moreover, there is an open subset $O\subset X$ such that all
points in $O$ verify that $W^s_\eps(x)=W^u_\eps(x)=\{x\}$
\par
Let us consider the following example, introduced in \cite{rr}:
Take a linear Anosov automorphism $f$ on ${\Bbb T}^2$, and choose
a fixed point $p\in {\Bbb T}^2$. Choose a stable and an unstable
separatrix of $p$, say $W^s_+(p)$ and $W^u_+(p)$. Lift copies
$\tilde p$, $\widetilde { W^s_+(p)}$, and $\widetilde{W^u_+(p)}$
of $p$, $W^s_+(p)$ and $W^u_+(p)$ respectively, in ${\Bbb R}^3$,
so that $\tilde{p}\in
\widetilde{W^s_+(p)}\cap\widetilde{W^u_+(p)}$ and
$\widetilde{W^\sigma_+(p)}$ be asymptotic to $W^\sigma_+(p)$,
$\sigma =s,u$. Extending $f$ continuously to this new set, gives
an expansive homeomorphism $F$. Observe that $F$ is an example of an
expansive homeomorphism on a continuum, having sinks and sources,
none of which is periodic. See Question \ref{pregunta.sink.en.cuenca de periodico}.
\par%
Take now any small perturbation of $F$ at $\tilde p$, say $G$, so
that there are no fixed points in
$O=\widetilde{W^s_+(p)}\cup\widetilde{W^u_+(p)}$, and so that
$\omega(x)\cup\alpha(x)\subset {\Bbb T}^2$ for all points in $O$.
$G$ is an expansive homeomorphism on a continuum with sensitive
dependence on initial conditions; however,
$W^s_\eps(x)=W^u_\eps(x)=\{x\}$ for all $x$ in $O$.

\end{document}